\documentclass[12pt,a4paper]{article}

\usepackage{amsfonts}
\usepackage{amssymb}
\usepackage{url}
\def\IC{\mathbb{C}}
\def\IR{\mathbb{R}}

\def\proof{\par\noindent\textbf{Proof.} }
\def\trace{\mbox{trace}\,}

\newtheorem{proposition}{Proposition}
\newtheorem{lemma}[proposition]{Lemma}
\newtheorem{theorem}[proposition]{Theorem}
\newtheorem{corollary}[proposition]{Corollary}
\newtheorem{remark}[proposition]{Remark}

\newtheorem{example}[proposition]{Example}

\begin{document}
\title{Norms and CB Norms of Jordan Elementary Operators}
\author{Richard M. Timoney}
\maketitle

\begin{abstract}
We establish lower bounds for norms and CB-norms of
elementary operators on $\mathcal{B}(H)$. Our main result concerns
the operator $T_{a,b}x = axb + bxa$ and we show $\left\|T_{a,b}
\right\| \geq \|a\| \|b\|$, proving a conjecture of
M. Mathieu.
We also
establish some other results and formulae for $\left\|T_{a,b}
\right\|_{cb}$ and $\left\|T_{a,b} \right\|$ for special cases.
\end{abstract}

Our results are related to a problem of M. Mathieu
\cite{MMBCan,MMBAus} asking whether $\left\|T_{a,b}\right\| \geq c \|a\|
\|b\|$ holds in general with $c = 1$. We prove this in
Theorem~\ref{mainthm} below.

In \cite{MMBAus} the inequality
is established for $c = 2/3$ and the best known result to date is $c =
2(\sqrt{2} -1)$ as shown in \cite{SZ,BCSZ,MagTurn}. There are simple examples
which show that $c$ cannot be greater than 1 in general and there are
results which prove the inequality with $c=1$ in special cases. The
case $a^* =a$ and $b^* = b$ is shown in \cite{Mag} where it is deduced
from $\left\|T_{a,b}\right\|_{cb} = \left\|T_{a,b}\right\|$ under
these hypotheses.

The equality of the the CB norm and the operator
norm of $T_{a,b}$ also holds if $a,b$ are commuting normal operators.
See section~\ref{secCommute} below for references.

A result for $c=1$ is shown in \cite{BB} under the assumption
that $\|a + z b\| \geq \|a\|$ for all $z \in \IC$. In more general
contexts similar results (with varying values of $c$) are shown in
\cite{CR,BCSZ}.

As this manuscript was being written we learned of another proof of
the main result (\cite{BBR}), using rather different methods.
Thanks are due to M. Mathieu for drawing our attention to this reference.

\textbf{Acknowledgement.} Part of this work was done during a visit by
the author to the University of Edinburgh in the autumn of 2002.
A significant impetus to the work arose from discussions with Bojan
Magajna and Aleksej Turn\v{s}ek during a visit to Ljubljana in March
2003 and the author is very grateful to them for that. Thanks also to
P. Legi\v{s}a for finding the reference to \cite{Horn} below.

\section{Preliminaries}

We call $T \colon \mathcal{B}(H) \to \mathcal{B}(H)$ an elementary
operator if $T$ has a representation 
\[
T(x) = \sum_{i=1}^\ell a_i x b_i
\]
with $a_i, b_i \in \mathcal{B}(H)$ for each $i$. We cite
\cite{AraMatBk} for an exposition of many of the known results on
(more general) elementary operators and for other concepts we cite a number of
treatises on operator spaces including \cite{EffrosRuan,Paulsen,ChristSinc}.
In particular we will use the completely bounded (or CB) norm
$\|T\|_{cb}$ of an elementary operator, the operator norm $\|T\|$ and
the estimate in terms of the Haagerup tensor product norm $\|T\| \leq
\|T\|_{cb} \leq \left\| \sum_{i=1}^\ell a_i \otimes b_i \right\|_h$.

We recall that the Haagerup norm of an element
$w \in \mathcal{B}(H) \otimes \mathcal{B}(H)$
(of the algebraic tensor product) is defined by
\[
\|w\|_h^2 = \inf \left\| \sum_{i=1}^k a_i a_i^* \right\|
\left\| \sum_{i=1}^k b_i^* b_i \right\|
\]
where the infimum is over all representations $w = \sum_{i=1}^k a_i
\otimes b_i$. Moreover this infimum is achieved with both $k$-tuples
$(a_1, a_2, \ldots, a_k)$ and $(b_1, b_2, \ldots, b_k)$ linearly
independent.

Throughout $H$ denotes a (complex) Hilbert space and $\mathcal{B}(H)$
the algebra of bounded linear operators on $H$. For $x$  in the class
of Hilbert-Schmidt operators on $H$ we denote the Hilbert-Schmidt norm
by $\|x\|_2$ (so that $\|x\|_2^2 = \trace x^*x$).

\section{Lower bounds}
\label{seclb}

\begin{lemma}
\label{lemsymm}
Given linearly independent $a, b \in \mathcal{B}(H)$, we can find
$c_1, c_2 \in \mathcal{B}(H)$, $\delta_1, \delta_2 > 0$ and $z \in \IC
\setminus \{0\}$ so that $a\otimes b  + b \otimes a = c_1 \otimes c_1
+ c_2 \otimes c_2$, $c_1 = (za + z^{-1}b)/\sqrt{2}$,
$c_2 = i(za - z^{-1}b)/\sqrt{2}$ and
\[
\|a\otimes b + b \otimes a\|_h = 
\| \delta_1 c_1c_1^* + \delta_2 c_2 c_2^*\| 
=
\| \delta_1^{-1} c_1^*c_1 + \delta_2^{-1} c_2^*c_2 \| .
\]
\end{lemma}

\proof We know from general facts cited above that the Haagerup norm
infimum
for $w = a\otimes b + b \otimes a$ is realised via a representation $w
= a_1 \otimes b_1 + a_2 \otimes b_2$. Moreover, by scaling $a_i$ to
$\lambda a_i$ and $b_i$ to $\lambda^{-1}b_i$ for a suitable $\lambda$
we can arrange that
\[
\|w\|_h = \| a_1 a_1^* + a_2 a_2^* \| = \|b_1^* b_1 + b_2^* b_2\|.
\]
We adopt a convenient matrix notation
\[
w = [a , b] \odot [b, a]^t = [a_1, a_2] \odot [b_1, b_2]^t
\]
for the two tensor product expressions above ($t$ for transpose)
and note that all possible (linearly independent)
representations of $w$ take the form
\[
w  = [a_1', a_2'] \odot [b_1', b_2']^t
= ([a_1, a_2] \alpha) \odot (\alpha^{-1} [b_1, b_2]^t)
\]
for a $2 \times 2$ invertible scalar matrix $\alpha$. We use the
transpose notation also for the linear operation on the tensor product
that sends $a_1 \otimes b_1$ to $b_1 \otimes a_1$. Then we have
\[
w = w^t = [b_1, b_2] \odot  [a_1, a_2]^t = ([a_1, a_2] \alpha) \odot
([b_1, b_2](\alpha^{-1}) ^t)^t.
\]
From $[b_1, b_2] = [a_1, a_2] \alpha$ and $[a_1, a_2] \alpha^t = [b_1,
b_2]$ together with linear independence we get $\alpha = \alpha^t$
symmetric.

We can now express $\alpha = u \Delta u^t$ where $u$ is a unitary
matrix and $\Delta$
is a
diagonal matrix with positive diagonal entries $\delta_1^{-1},
\delta_2^{-1}$ (\cite[Takagi's factorisation, 4.4.4]{Horn} --- see
also the problems on pages 212, 217 in \cite{Horn}).
Take $[a_1', a_2'] = [a_1, a_2] u $,
$[b_1', b_2'] = [b_1, b_2] (u^{-1})^t$ so that
\[
w = [a_1', a_2'] \odot [b_1', b_2']^t,
\]
\[
\|w\|_h = \| (a_1')(a_1')^* + (a_2')(a_2')^* \|
= \| (b_1')^*(b_1') + (b_2')^*(b_2') \|
\]
and
\[
[a_1', a_2'] \Delta = [a_1, a_2] u \Delta = [a_1, a_2] \alpha
(u^{-1})^t = [b_1, b_2] (u^{-1})^t = [b_1', b_2'].
\]
In other words, $a_i' \delta_i^{-1} = b_i'$ ($i=1,2$).

We now take $c_i = \sqrt{\delta_i} b_i'$ and we then have $w = c_1 \otimes
c_1 + c_2 \otimes c_2$ together with
\[
\|w\|_h = \| \delta_1 c_1c_1^* + \delta_2 c_2c_2^* \|
= \| \delta_1^{-1} c_1^*c_1 + \delta_2^{-1} c_2^*c_2 \|.
\]
It remains to relate $c_1, c_2$ to $a, b$ as claimed. If we put $a' =
(c_1 - i c_2)/\sqrt{2}$ and $b' = (c_1 +ic_2)/\sqrt{2}$ we have
\[
w = a' \otimes b' + b' \otimes a' = [a', b'] \odot [b',a']^t = [a,b]
\odot [b,a]^t.
\]
An easy argument shows that there is $z \in \IC$ with either
$a' = za$ and $b' = z^{-1}b$ or else $a' = z^{-1}b$ and $b' = za$.
The first case is exactly as required but for the second
case we need to swap the roles of $c_1$ and $c_2$.

\begin{theorem}
\label{L2bound}
Assume that $H$ is two-dimensional and $a,b \in \mathcal{B}(H)$.

Let $T_{a,b}(x) = axb + bxa$. Then
\[
\left\|T_{a,b}\right\|_{cb} \geq \|a\|_2 \|b\|_2.
\]
\end{theorem}

\proof In the case where $a,b$ are linearly dependent ($a = \lambda
b$, say, $T_{a,b}x = 2 \lambda a x a$) we know $\|T\|_{cb} = \|T\| =
2\|a\| \|b\| \geq \|a\|_2 \|b\|_2$. So we deal only with the case of
independent $a,b$.

We first apply Lemma~\ref{lemsymm}, $\left\|T_{a,b}\right\|_{cb} = \|a
\otimes b + b \otimes a\|_h$ and the fact that the norm of a $2
\times 2$ positive matrix (the max of the eigenvalues) is at least half
the trace to get
\begin{eqnarray*}
\left\|T_{a,b}\right\|_{cb} &\geq &
\frac{1}{2} \left( \delta_1 \|c_1\|_2^2 + \delta_2 \|c_2\|_2^2 \right)\\
\left\|T_{a,b}\right\|_{cb} &\geq &
\frac{1}{2} \left( \delta_1^{-1} \|c_1\|_2^2 + \delta_2 ^{-1}\|c_2\|_2^2
\right)\\
\end{eqnarray*}
We deduce
\begin{eqnarray*}
\left\|T_{a,b}\right\|_{cb} &\geq &
\frac{1}{4}  \left( (\delta_1 + \delta_1^{-1}) \|c_1\|_2^2 +
(\delta_2 +  \delta_2 ^{-1}) \|c_2\|_2^2 \right)\\
&\geq & \frac{1}{2} \left( \|c_1\|_2^2 +  \|c_2\|_2^2 \right)\\
&=& \frac{1}{2} \trace \left( c_1^* c_1 + c_2^* c_2 
\right)\\
&=& \frac{1}{2} \trace \left( (za)^* (za) + 
(z^{-1}b)^* (z^{-1}b) \right) \\
&=& \frac{1}{2} \left( \|za\|_2^2 + \|z^{-1}b\|_2^2 \right)\\
&\geq&  \|za\|_2 \|z^{-1}b\|_2 = \|a\|_2 \|b\|_2.
\end{eqnarray*}

\begin{corollary}[\cite{MagTurn}, Theorem 2.1] For $a, b \in
\mathcal{B}(H)$ ($H$ arbitrary)
\[
\left\| T_{a,b} \right\|_{cb} \geq \|a\| \|b\|.
\]
\end{corollary}

\proof
We can reduce the proof to the case where $H$ is two-dimensional by
the argument given in \cite[Theorem 2.1]{MagTurn} (take unit vectors
$\xi, \eta \in H$ where $\|a \xi\| \geq \|a \|- \varepsilon$ and
$\|b\eta\| \geq \|b\|- \varepsilon$; consider $T_{qap, qbp}$ where $p$
is a projection onto the span of $\xi, \eta$ and $q$ a projection onto
the span of $a\xi, b\eta$). In two dimensions the result follows from
Theorem~\ref{L2bound}.

\begin{proposition}
\label{symmprop}
If $a, b \in \mathcal{B}(\IC^2)$ are symmetric matrices, then
\[
\left\| T_{a,b} \right\|_{cb} =
\left\| T_{a,b} \right\|
= \inf_{x > 0} \left \| x a a^* + (1/x) b b^* \right\|
\]
\end{proposition}

\proof
Now $c_1, c_2$ obtained from Lemma~\ref{lemsymm} are symmetric
matrices.
Using $c_i^* = \bar{c}_i =$ the complex conjugate matrix we have
\[
\left\| \delta_1^{-1} c_1^*c_1 + \delta_2^{-1} c_2^*c_2 \right\| 
= \left\| \delta_1^{-1} \bar{c}_1c_1 + \delta_2^{-1} \bar{c}_2c_2 \right\|
= \left\| \delta_1^{-1} c_1\bar{c}_1 + \delta_2^{-1} c_1\bar{c}_2 \right\|
\]
Thus
\begin{eqnarray*}
\left\| T_{a,b} \right\|_{cb} &\geq&
\left\| \frac{\delta_1 +\delta_1^{-1}}{2} c_1c_1^* + \frac{\delta_2 +
\delta_2^{-1}}{2} c_2 c_2^*\right\| \\
&\geq &
\|  c_1c_1^* + c_2 c_2^*\| \\
&=& \|  c_1^*c_1 + c_2^* c_2\|
\end{eqnarray*}
so that the infimum in the Haagerup tensor norm is attained with
$\delta_1 = \delta_2 = 1$. We thus have
\[
\left\| T_{a,b} \right\|_{cb} = \inf_{z} \left\|
|z|^2 a a^* + |z|^{-2} b b^* \right\|
\]
and the desired formula for $\left\| T_{a,b} \right\|_{cb}$ (taking $x
= |z|^2$).

From \cite{TimoneyNE} we know that the
convex hulls of the following two sets of matrices intersect
\begin{eqnarray}
W_l &=& \left\{  \left[ \begin{array}{cc}
\langle c_1 c_1^* \xi, \xi \rangle &
\langle c_2 c_1^* \xi, \xi \rangle \\
\langle c_1 c_2^* \xi, \xi \rangle &
\langle c_2 c_2^* \xi, \xi \rangle 
\end{array} \right] : \xi \in H, \|\xi\| = 1, \right. \nonumber\\
&&
\label{eqnNRl}
\qquad \qquad \left. \left\langle
\left(\sum_{i=1}^2 c_i c_i^*\right) \xi, \xi \right \rangle =
\left\| T_{a,b} \right\|_{cb} \right\},
 \\
W_r &=& \left\{  \left[ \begin{array}{cc}
 \langle c_1^* c_1 \eta, \eta \rangle &
\langle c_2^* c_1 \eta, \eta \rangle \\
 \langle c_1^* c_2 \eta, \eta \rangle &
 \langle c_2^* c_2 \eta, \eta \rangle 
\end{array} \right] : \eta \in H, \|\eta\| = 1, \right.
 \nonumber\\
&&
\label{eqnNRr}
\qquad \qquad \left. \left\langle
\left(\sum_{i=1}^2 c_i^* c_i\right) \eta, \eta \right \rangle =
\left\| T_{a,b} \right\|_{cb} \right\}.
 \end{eqnarray}
Moreover the equality $\left\| T_{a,b} \right\|_{cb} =
\left\| T_{a,b} \right\|$ holds if and only if the sets themselves
intersect. For either of the sets (say $W_l$) to consist of more than
one element, the hermitian operator concerned must have a double
eigenvalue of the maximum eigenvalue $\left\| T_{a,b} \right\|_{cb}$,
which means that (taking the case $W_l$)
\[
\sum_{i=1}^2 c_i c_i^*
\]
is a multiple of the $2 \times 2$ identity matrix. But then by complex
conjugation and symmetry $\sum_{i=1}^2 c_i^* c_i $ is the same
multiple of the identity.

In the case when $W_l$ (and $W_r$ by the symmetry) are singletons, we
have $\left\| T_{a,b} \right\|_{cb} =
\left\| T_{a,b} \right\|$ and
using the following lemma, we can complete the proof for the other
case.

\begin{lemma}
\label{lemAdoublesided}
If $c_1, c_2 \in \mathcal{B}(\IC^2)$ are symmetric and
satisfy
$c_1c_1^* + c_2c_2^* 
= $ a multiple of
the identity matrix,
there exists $u$ unitary so that either
$uc_1u^t$ and $u c_2 u^t$ are both diagonal ($t$ for transpose) or
\[
uc_1u^t = 
\left( \begin{array}{cc}
\lambda & 0\\
0 & \lambda \end{array} \right), \quad
u c_2 u^t =
\left( \begin{array}{cc}
\zeta \alpha & \zeta \beta\\
\zeta \beta & -\zeta \bar{\alpha} \end{array} \right)
\]
with $\lambda > 0$, $\beta > 0$, $|\zeta| = 1$.
\end{lemma}

\proof
We can find $u$ so that $uc_1u^t$ is diagonal (with positive entries,
\cite[4.4.4]{Horn}).

We can replace $c_i$ by $uc_iu^t$ ($i=1,2$) and assume without loss of
generality that $c_1$ is diagonal. Then $c_2 c_2^*$
is diagonal, which means that the rows of $c_2$ are
orthogonal. An easy analysis shows that either $c_2$ is diagonal
or is a multiple (of modulus one) of a matrix of the form
\[
\left( \begin{array}{cc}
\alpha & \beta\\
\beta & -\bar{\alpha} \end{array} \right)
\]
The relation satisfied by $c_1$ and
$c_2$ dictates that $c_1$ is a multiple of the identity in the latter
case.

\proof (of Proposition~\ref{symmprop}, completed).
Invoking the lemma and the fact that $S(x) = u T(u^t x u) u^t$ has
the same norm as $T$, and the same CB norm, we can reduce to the case
where $c_1, c_2$ generate a commutative $C^*$ algebra. In this case
the fact that $\|S\|_{cb} = \|S\|$ is known (see references in
section~\ref{secCommute}).

\begin{theorem}
\label{mainthm}
If $a,b \in \mathcal{B}(H)$ and
$T_{a,b}(x) = axb + bxa$. Then
\[
\left\|T_{a,b}\right\| \geq \|a\| \|b\|.
\]
More generally, the same inequality holds 
if $A$ is a prime $C^*$-algebra, $a,b$ are in
the multiplier algebra of $A$ and $T_{a,b} \colon A \to A$ is
$T_{a,b}(x) = axb + bxa$.
\end{theorem}

\proof
As shown in \cite{MMBAus} and \cite[Theorem 2.1]{MagTurn}, the
essential case is the case where $A = \mathcal{B}(H)$ and
$H = \IC^2$ is 2-dimensional. We show
in this case that $\left\|T_{a,b}\right\| \geq \|a\| \|b\|_2 \geq \|a\|
\|b\|$ and so we
can assume $\|a\| = \|b\|_2 = 1$ ($a, b \in \mathcal{B}(\IC^2)$).

There exists $u, v$ unitary so that $uav$ is a diagonal matrix with
diagonal entries $1, \lambda$, $0 \leq |\lambda| \leq 1$.
Replacing $T$ by $S(x) = u T(vxu)v$ we can assume that
\[
a =
\left( \begin{array}{cc}
1 & 0\\
0 & \lambda \end{array} \right), \quad
b = \left( \begin{array}{cc}
b_{11} & b_{12}\\
b_{21} & b_{22} \end{array} \right).
\]
By multiplying $b$ by a scalar of modulus 1 we can assume that $b_{12}
= |b_{12}|$. Multiplying both $a$ and $b$ by a diagonal unitary $u$ with 
diagonal entries 1 and $\bar{b}_{21}/|b_{21}|$  (that is,
replacing $T$ by $S(x) = u T(xu)$) we can assume also that
$b_{21} = |b_{21}|$.

Now consider $T_t(x) = T(x^t)^t = a x b^t + b^t x a$ and
\[
T_s(x) = \frac{1}{2} \left( T(x) + T_t(x) \right) = a x b_s + b_s x
a
\]
with
\[
b_s = \frac{1}{2}(b + b^t) 
= \left( \begin{array}{cc}
b_{11} & s_{12}\\
s_{12} & b_{22} \end{array} \right), \quad s_{12} = \frac{b_{12} +
b_{21}}{2}.
\]

We claim that $\|T_s\| \geq 1$ and this will prove the theorem because
$\|T_t\| = \|T\|$ and so $\|T_s\| \leq \|T\|$.

To show $\|T_s\| \geq 1$ we invoke Proposition~\ref{symmprop} and show
$\|T_s\|_{cb} \geq 1$. Note
\[
\frac{1}{2} \leq \|b_s\|_2^2 = \|b\|_2^2 -
\frac{1}{2}(b_{12} - b_{21})^2 \leq 1,
\]

\[
b_s b_s^* = \left( \begin{array}{cc}
|b_{11}|^2 + s_{12}^2 & s_{12}(b_{11} + \bar{b}_{22}) \\
s_{12}(\bar{b}_{11} + b_{22}) & |b_{22}|^2 +
s_{12}^2
\end{array} \right)
\]
and write
$\mu_i^2 = |b_{ii}|^2 + s_{12}^2$ ($i = 1,2$) for the diagonal entries.

Now consider a unit vector $\xi = (\xi_1, \xi_2) \in \mathbb{C}^2$.
Then
\begin{eqnarray*}
\| x a a^* + (1/x) b_s b_s^* \|
&\geq& \langle (x a a^* + (1/x) b_s b_s^*) \xi, \xi \rangle\\
&=& x \langle a a^* \xi, \xi \rangle + (1/x)  \langle b_s b_s^* \xi, \xi \rangle\\
&\geq& 2 \sqrt{ \langle a a^* \xi, \xi \rangle \langle b_s b_s^* \xi, \xi
\rangle}
\end{eqnarray*}
and we claim that there is a point in the joint numerical
range
\[
W = \{ (x, y) = ( \langle a a^* \xi, \xi \rangle, \langle b_s b_s^* \xi, \xi 
\rangle ) : \|\xi\| = 1 \} \subseteq \IR^2
\]
which is also on (or above) the hyperbola $xy = 1/4$. Verifying the claim will
complete the proof.

\emph{We assume  from now on} that $\lambda = 0$,
as this is the hardest case (smallest
$\langle a a^* \xi, \xi \rangle$).

Being the joint numerical range of two hermitian operators (or the
numerical range of the single operator $a a^* + i  b_s b_s^*$), $W$ is
a convex set in the plane. In fact, because the space is
2-dimensional, $W$ is either a straight line (in
the case where the two operators commute, that is $s_{12}(b_{11} +
\bar{b}_{22}) = 0$) or else an ellipse (together with its interior)
\cite[I.6.2]{Bhatia}.
The ellipse touches the
vertical lines $x=0$ and $x=1$ at the points $(0, \mu_2^2)$ and $(1,
\mu_1^2)$. Hence the centre of the ellipse is at the midpoint $(x_0,
y_0) = (1/2, (1/2)(\mu_1^2 + \mu_2^2)) = (1/2, (1/2) (|b_{11}|^2 +
|b_{22}|^2) + s_{12}^2) = (1/2, (1/2) \|b_s\|_2^2)$.

In the case where we have a line and not a genuine ellipse, either
$s_{12} = 0$ (then the midpoint is $(1/2,
1/2)$ and so on the hyperbola) or $b_{11} = -\bar{b}_{22}$ and the
line is horizontal (at $y = (1/2) \|b_s\|_2^2 \geq 1/4$ and so also
meets the hyperbola).
If $|b_{11}| \geq |b_{22}|$, then the point $(x,y) =
(1, \mu_1^2)$ on the
ellipse already satisfies $4xy \geq 1$ and so we assume that $|b_{22}|
> |b_{11}|$.

For the genuine ellipse case we write its equation in the form
\begin{equation}
\label{ellipseeqn}
\alpha_{11} (x - x_0)^2 + 2 \alpha_{12} (x - x_0)(y-y_0) + (y - y_0)^2 +
\beta = 0.
\end{equation}
Using the information that the ellipse has a vertical tangent at $(0,
\mu_2^2)$ and its intersection with the line $x = 1/2$ is the line
segment
$\{(1/2, y) :|y-y_0| \leq s_{12} |b_{11} +
\bar{b}_{22}| \}$ (take $\xi$ with $\xi_1 = 1/\sqrt{2}$),
we can solve for the coefficients
\begin{eqnarray}
\alpha_{12} &=& \mu_2^2 - \mu_1^2 = |b_{22}|^2 - |b_{11}|^2 \nonumber\\
\label{betaformula}
\beta &=& - s_{12}^2 |b_{11} + \bar{b}_{22}|^2\\
\alpha_{11} & =&  (|b_{11}|^2 - |b_{22}|^2)^2 + 4 s_{12}^2
|b_{11} + \bar{b}_{22}|^2 = \alpha_{12}^2 - 4 \beta \nonumber
\end{eqnarray}

We can rewrite the equation in the form
\[
(\alpha_{12} (x - x_0) + (y - y_0))^2 - 4 \beta (x - x_0)^2 + \beta =
0
\]
and so we can parametrise the ellipse via
\begin{eqnarray}
\label{xgeqn}
x & = & x_0 + (1/2) \sin \omega \\
y & = & y_0 - (1/2) \alpha_{12} \sin \omega + \sqrt{- \beta} \cos \omega
\nonumber\\
&=& (1/2) (|b_{11}|^2 + |b_{22}|^2) + s_{12}^2 - (1/2) (|b_{22}|^2 -
|b_{11}|^2) \sin \omega \nonumber\\
&& \quad + s_{12} |b_{11} + \bar{b}_{22}| \cos \omega
\label{ygeqn}
\end{eqnarray}
($0 \leq \omega \leq 2 \pi$.) We look for $\omega \in [0, \pi/2]$
where $4xy \geq 1$. We use $|b_{11} + \bar{b}_{22}| \geq |b_{22}| -
|b_{11}| = \epsilon_{12}$ (say) and represent for convenience
$|b_{11}|^2 + |b_{22}|^2 = \cos^2 \theta$ ($0 \leq \theta < \pi/2$).
Note
$4s_{12}^2 \geq (b_{12} - b_{21})^2$, $2 s_{12}^2 \geq
(1/2) (b_{12} - b_{21})^2 = 1 - \| b_s\|_2^2$, $4 s_{12}^2 \geq 1 -
\cos^2 \theta$ and $s_{12} \geq (1/2) \sin \theta$. Moreover
$|b_{22}| + |b_{11}| \leq \sqrt{2} \cos \theta$. Thus
\begin{equation}
\label{twoylb}
2 y \geq (1/2) + (1/2) \cos^2 \theta + \epsilon_{12}
( \sin \theta \cos \omega - \sqrt{2} \cos \theta \sin \omega )
\end{equation}
Choose $\omega = \tan^{-1} ( (1/\sqrt{2}) \tan \theta)$, $\sin \omega
= \sin \theta /\sqrt{\sin^2 \theta + 2 \cos^2 \theta}$ and
\[
4xy \geq \left ( 1 + \frac{\sin \theta }{\sqrt{1 + \cos^2 \theta}} \right)
(1/2 + (1/2) \cos^2 \theta) \geq 1.
\]

\begin{remark}
With some additional effort, we can adapt the proof above to establish
the lower bound $\left\|T_{a,b}\right\| \geq \|a\|_2 \|b\|_2$ for the
case $a, b \in \mathcal{B}(\IC^2)$ (and thus get a stronger result
than Theorem~\ref{L2bound}).

\rm
It seems that this does not follow from the methods used in \cite{BBR}.

\proof
A sketch of the additional details follows. We assume by symmetry
that $\|a\|_2/\|a\|
\leq \|b\|_2/\|b\|$ and normalise $\|a\| =1$, $\|b\|_2 = 1$ as before.
This time we cannot assume $\lambda = 0$, but we note that $|\det b|
\geq |\lambda|/(1 + |\lambda|^2)$ (for example, take $b = u b_0 v$
where $u,v$ are unitary and $b_0$ is diagonal with diagonal entries
$1/\sqrt{1 + \mu^2}$ and $\mu/\sqrt{1 + \mu^2}$, $1 \geq \mu \geq |\lambda|$).

In this case the ellipse will have vertical tangents at $x =
|\lambda|^2$ and $x= 1$ and will be centered at $(x_0, y_0) = ((1 +
|\lambda|^2)/2, (1/2) \|b_s\|_2^2)$. The equation (\ref{ellipseeqn})
of the ellipse now has
\[
\alpha_{12} = \frac{|b_{22}|^2 - |b_{11}|^2}{1 - |\lambda|^2},
\]
$\beta$ as in (\ref{betaformula}) and
$\alpha_{11} = \alpha_{12}^2 - 4 \beta/(1 - |\lambda|^2)^2$. We can rewrite the equation
of the ellipse as
\[
(\alpha_{12} (x - x_0) + (y - y_0))^2 - \frac{4 \beta}{(1 -
|\lambda|^2)^2} (x - x_0)^2 + \beta = 0
\]
and then we can parametrise via 
\begin{equation}
\label{newxg}
x  =  (1/2)(1 + |\lambda|^2) + (1/2) (1 - |\lambda|^2) \sin \omega 
\end{equation}
(in place of (\ref{xgeqn})) and (\ref{ygeqn}) as before.

We now seek a point $(x,y)$ on the ellipse where $4xy \geq 1 +
|\lambda|^2$.

To dispose of the case $|b_{11}| \geq |b_{22}|$ we show $4 y_0 \geq  1
+ |\lambda|^2$ (and this also deals with the case where the ellipse
degenerates into a line). Using $\|b\|_2 = 1$,
\begin{eqnarray*}
4 y_0 &=& 2 \|b_s\|^2 = 2 - (b_{12} - b_{21})^2 = 1 + (|b_{11}|^2 +
|b_{22}|^2 + 2 b_{12} b_{21})
\\
&\geq& 1 + 2 | b_{11} b_{22} - b_{12}
b_{21}| \geq 1 + 2 \frac{|\lambda|}{1 + |\lambda|^2} \geq 
1 + |\lambda|^2.
\end{eqnarray*}

When $\epsilon_{12} = |b_{22}| -
|b_{11}| > 0$
we choose the same $\omega$ as before.
From the
lower bound (\ref{twoylb}) and (\ref{newxg}) 
we get the desired $4xy \geq 1 +
|\lambda|^2$ if we have $ \cos^2 \theta \geq 2 |\lambda|^2/(1 +
|\lambda|^4)$. For the remaining case note that
\[
2 y \geq |b_{11}|^2 + |b_{22}|^2 + 2 s_{12}^2
= \frac{1}{2} + \frac{1}{2} (|b_{11}|^2 + |b_{22}|^2) + b_{12} b_{21}
\geq \frac{1}{2} + |\det b|
\]
and the resulting $2y \geq 1/2 + |\lambda|/(1 + |\lambda|^2)$ is a
better lower bound that (\ref{twoylb})
when $\cos^2 \theta < 2 |\lambda|/(1 +
|\lambda|^2)$. In this situation we do get $4xy \geq  1 +
|\lambda|^2$. All eventualities are now covered because $2
|\lambda|^2/(1 + |\lambda|^4) \leq 2 |\lambda|/(1 + |\lambda|^2)$.
\end{remark}

\section{Commuting cases}
\label{secCommute}

We consider now some cases where we can find relatively explicit
formulae for $\|T_{a,b}\|$. These may shed some light on the
difficulty of finding any explicit formula for the norm of a general
elementary operator. One may consider the Haagerup formula for the
CB norm as an explicit formula, though we shall observe that this is
not so simple to compute even in the simplest cases.

The equality of the CB norm and the operator
norm of $T_{a,b}$ holds if $a,b$ are commuting normal operators.
This appears already in the unpublished \cite{Haagerup}. A significant
part of the argument from \cite{Haagerup} is published in
\cite[\S 5.4]{AraMatBk} and the remaining part uses the fact that all
states on a commutative $C$*-algebra are vector states. (By the
Putnam-Fuglede theorem the $C$*-algebra generated by commuting normal
operators is commutative.) See also \cite[Theorem 2.1]{Smith} for a
more general result on bimodule homomorphisms.
Another proof (with slightly weaker hypotheses)
is in \cite{TimoneyNE}.

We deal here only with $H$ of dimension 2.

\begin{proposition}
\label{propCommute}
If $H$ is two-dimensional and $a,b \in
\mathcal{B}(H)$ commute, then $\left\| T_{a,b} \right\|_{cb} = 
\left\| T_{a,b} \right\|$.
\end{proposition}

\proof
We can find an orthonormal basis of $H$ so that $a$ and
$b$ both have upper triangular ($2 \times 2$) matrices. If
$a,b$ are diagonal, then they generate a commutative $C$*-subalgebra
of $\mathcal{B}(H)$ and in this case that $\left\|
T_{a,b} \right\|_{cb} = \| a \otimes b + b \otimes a \|_h = \left\|
T_{a,b} \right\|$ (see above).

Now $c_1, c_2$ obtained from Lemma~\ref{lemsymm} are also commuting
upper triangular matrices. As used already in
(\ref{eqnNRl}) -- (\ref{eqnNRr}), from \cite{TimoneyNE} we know that the
convex hulls of the two sets of matrices intersect. In this case the
sets are as not quite as before. Each $c_i$ should be
replaced by $\sqrt{\delta_i} c_i$ in
the definition of $W_l$ and by $1/\sqrt{\delta_i} c_i$ for $W_r$.
Moreover the equality $\left\| T_{a,b} \right\|_{cb} =
\left\| T_{a,b} \right\|$ holds if and only if the sets themselves
intersect. For either of the sets (say $W_l$) to consist of more than
one element, the hermitian operator concerned must have a double
eigenvalue of the maximum eigenvalue $\left\| T_{a,b} \right\|_{cb}$,
which means that (taking the case $W_l$)
\[
\sum_{i=1}^2 \delta_i c_i c_i^*
\]
is a multiple of the $2 \times 2$ identity matrix. But the following
lemma asserts that this cannot happen unless $\sqrt{\delta_1}c_1$ and
$\sqrt{\delta_2}c_2$ are simultaneously diagonalisable (the case
where we know the result). So $W_l$ and $W_r$ have one element each,
they intersect and the result follows.

\begin{lemma} If $a_1, a_2$ are commuting elements of $\mathcal{B}(H)$
with $H$ of dimension $2$ and if $a_1a_1^*+a_2a_2^*$ is a multiple of
the identity, then $a_1, a_2$ generate a commutative *-subalgebra of
$\mathcal{B}(H)$.
\end{lemma}

\proof In a suitable orthonormal basis for $H$ we can represent $a_1,
a_2$ as upper triangular matrices
\[
a_1 = \left[ \begin{array}{cc} x_1 & y_1 \\ 0 & z_1 \end{array}
\right],
\qquad
a_2 = 
\left[ \begin{array}{cc} x_2 & y_2 \\ 0 & z_2 \end{array} 
\right]
\]
and then the condition for them to commute is  $y_1(x_2 - z_2) =
y_2(x_1 - z_1)$. (For later reference we call this value $\rho$).
So if $y_1 = 0$, then either $y_2$ also zero (both
matrices diagonal and we are done) or else $x_1 = z_1$ and $a_1 = x_1
I_2$ is a multiple of the identity. But then $a_2a_2^*$ is a multiple
of the identity and this forces $y_2 = 0$ (both diagonal again).

In the case when $y_1$ and $y_2$ are both nonzero, we compute
\[
a_1 a_1^* + a_2 a_2^* = \left[ \begin{array}{cc}
|x_1|^2 + |y_1|^2 + |x_2|^2 + |y_2|^2 & y_1 \bar{z}_1 + y_2
\bar{z}_2\\
\bar{y}_1 z_1 + \bar{y}_2 z_2 & |z_1|^2 + |z_2|^2
\end{array} \right]
\]
Thus we have $y_1 \bar{z}_1 + y_2 \bar{z}_2 = 0$, which implies
$(z_1, z_2) = \omega (\bar{y}_2, -\bar{y}_1)$ for some $\omega \in
\IC$. We also have equality of the two diagonal entries of the above
matrix which gives us
\[
|x_1|^2 + |x_2|^2 = (|\omega|^2 -1) (|y_1|^2 + |y_2|^2)
\]
Now $x_1 = \rho/y_2 + z_1 = \rho/y_2 + \omega \bar{y}_2$ and $x_2 =
\rho/y_1 - \omega \bar{y}_1$, yielding
\[
\left| \frac{\rho}{y_2} + \omega \bar{y}_2 \right|^2
+\left| \frac{\rho}{y_1} - \omega \bar{y}_1 \right|^2
= (|\omega|^2 -1)(|y_1|^2 + |y_2|^2).
\]
and hence the impossible condition
\[
|\rho|^2 (|y_1|^{-2} + |y_2|^{-2}) = -(|y_1|^2 + |y_2|^2)
\]

\begin{example}
\rm
Consider $T_{a,b}$ acting on $\mathcal{B}(\IC^2)$
with $a,b$ diagonal $2 \times 2$ matrices. Then $c_1$, $c_2$ in Lemma
\ref{lemsymm} are also diagonal and we can see then directly that
\[
\|c_1 c_1^* + c_2 c_2^*\| \leq \frac{1}{2} ( \| \delta_1 c_1c_1^* +
\delta_2 c_2c_2^* \| + \| \delta_1^{-1} c_1^*c_1 + \delta_2^{-1} c_2^*
c_2 \| )
\]
so that the Haagerup norm is minimised with $\delta_1 = \delta_2 = 1$.
Also $\|c_1 c_1^* + c_2 c_2^*\| = \| |z|^2aa^* + |z|^{-2} b b^* \|$
and so the Haagerup norm is the minimum of this.

Say the diagonal entries are $\lambda_1, \lambda_2$ for $a$ and
$\mu_1, \mu_2$ for $b$. Normalising $a$ and $b$ to have norm one, we
can assume $\max(|\lambda_1|,|\lambda_2|) = 1$ and $\max(|\mu_1|,
|\mu_2|) = 1$. If they both attain the maximum at the same index then
it is easy to see that $\left\| T_{a,b} \right\| = 2 = 2 \|a\| \|b\|$.
If not, assume by symmetry that $|\lambda_1| = 1 = |\mu_2|$ and that
$|\mu_1| \leq |\lambda_2|$.
The Haagerup norm is then the minimum value of the maximum of two
functions, and can be computed by elementary means. It
gives the norm (the same as the CB norm in this case) as
\begin{equation}
\left\|T_{a,b} \right\|
=
\left\{ \begin{array}{ll}
2 |\lambda_2| & \mbox{if } |\lambda_2|  \geq 1/\sqrt{2}\\
&
\mbox{ and } |\mu_1|^2 < 2 - |\lambda_2|^{-2}\\
\displaystyle
\frac{1 - |\mu_1|^2 |\lambda_2|^2}{
\sqrt{(1 - |\mu_1|^2)(1 - |\lambda_2|^2)} }& \mbox{otherwise}
\end{array} \right.
\end{equation}

Summarising the calculation in a basis independent way, we can state
the following.
\end{example}

\begin{proposition}
Suppose 
that $a,b \in \mathcal{B}(\IC^2)$ are commuting normal operators
and that $\|a\|_2/\|a\| \geq \|b\|_2/\|b\|$. If $a$, $b$ attain their
norms at a
common unit vector,
then $\left\|T_{a,b} \right\| = 2\|a\| \|b\|$. If not
\begin{equation}
\left\|T_{a,b} \right\| =
\left\{ \begin{array}{ll}
2 \|b\| \sqrt{\|a\|_2^2 -\|a\|^2} \\
\qquad \mbox{ if } \|a\|_2 \geq  \sqrt{3/2} \|a\|
\\
 \qquad \mbox{ and } \|b\|_2^2 < 3\|b\|^2 - (\|a\|^2\|b\|^2)/(\|a\|_2^2 -\|a\|^2)\\
\displaystyle
\frac{\|a\|_2^2\|b\|^2 + \|a\|^2\|b\|_2^2 - \|a\|_2^2 \|b\|_2^2 }
{\sqrt{(2\|a\|^2- \|a\|_2^2)(2\|b\|^2-\|b\|_2^2)}}
\qquad
\mbox{otherwise}
\end{array} \right.
\end{equation}
\end{proposition}

\proof
Note that in a suitable orthonormal basis of $\IC^2$, $a,b$ will both
be represented by diagonal matrices.

\section{A formula for self-adjoint operators}

Our aim here is to present a proof of a formula from \cite{Mag} that
follows a similar approach to the one used in
section~\ref{seclb}.

For a linear operator $T \colon \mathcal{B}(H) \to \mathcal{B}(H)$
we denote by $T^*$ the associated operator defined by $T^*(x) =
T(x^*)^*$. We call $T$ self-adjoint if $T^* = T$.

\begin{lemma}[\cite{TimoneyNE}]
\label{lemmStarSymm}
For $T \colon \mathcal{B}(H) \to \mathcal{B}(H)$ a self-adjoint elementary
operator, there is a representation $Tx = \sum_{i=0}^\ell
\varepsilon_i c_i x c_i^*$ with $c_i \in \mathcal{B}(H)$,
$\varepsilon_i \in \{ -1, 1\}$ for each $i$ and
\[
\|T\|_{cb} = \left\| \sum_{i=1}^\ell c_i c_i^*\right\|.
\]
\end{lemma}

\begin{lemma}[\cite{TimoneyNE}]
\label{lemmSignature}
Let $T = T^* \colon \mathcal{B}(H) \to \mathcal{B}(H)$ be an
elementary operator, $Tx = \sum_{i=1}^k
c_ixc_i^* - \sum_{i=k+1}^\ell c_ixc_i^*$ with $0 \leq k \leq \ell$
and $(c_i)_{i=1}^\ell$ linearly independent. (We include $k=0$ for the
case where the first summand is absent and when $k = \ell$ the second
summand is absent.) Then the ordered pair $(k, \ell-k)$ (which we
could call the `signature') is the same for all such representations
of $T$.
\end{lemma}

\begin{example}[\cite{Mag}]
For $T \colon \mathcal{B}(H) \to \mathcal{B}(H)$
given by $Tx = axb^* + bxa^*$ with
$a$, $b$ linearly independent, we have
\[
\|T\|_{cb} = \inf \left\{ \| raa^* + sbb^*  +2t \Im(ab^*)\| : r >
0,s>0, t \in \IR, rs -t^2 =1 \right\}
\]
(where $\Im(ab^*) = (ab^* - ba^*)/(2i)$ is the imaginary part).
\end{example}

\proof We can rewrite $Tx = c_1 x c_1^* - c_2xc_2^*$ if we take
$c_1 = (a+b)/\sqrt{2}$ and $c_2 = (a-b)/\sqrt{2}$. Note for later use
that we can undo this change by $a = (c_1+c_2)/\sqrt{2}$, $b =
(c_1-c_2)/\sqrt{2}$.

According to Lemma~\ref{lemmStarSymm} and
Lemma~\ref{lemmSignature} we can find $\|T\|_{cb}$ as the infimum of
$\|c_1'(c_1')^* + c_2'(c_2')^*\|$ where
\[
[c_1', c_2'] =[c_1, c_2]\alpha
\]
and $\alpha$ is an invertible $2 \times 2$ matrix with the property
that
\[
\alpha
\left[
\begin{array}{cc}
1 & 0\\
0 & -1
\end{array} \right] \alpha^*
= \left[
\begin{array}{cc}
1 & 0\\
0 & -1
\end{array} \right].
\]
As unitary diagonal $\alpha$ have no effect on the estimate
$\|c_1'(c_1')^* + c_2'(c_2')^*\|$ we can work modulo these unitaries
and then
elementary analysis of the possibilities shows that we need only
consider the cases
\[
\alpha = \left[
\begin{array}{cc}
p & \sqrt{p^2-1} e^{i \theta}\\
\sqrt{p^2-1} e^{-i \theta} & p
\end{array} \right]
\]
(with $p \geq 1$, $\theta \in \IR$). This leads us to consider only
\[
[c_1', c_2'] = [pc_1 + \sqrt{p^2-1} e^{- i \theta} c_2,
\sqrt{p^2-1} e^{i \theta} c_1 + p c_2].
\]
Hence
\begin{eqnarray*}
\|T\|_{cb}
&=&\inf_{p \geq 1, \theta \in \IR} \|c_1'(c_1')^* + c_2'(c_2')^*\|\\
&=& \inf \| (2p^2 -1) (c_1c_1^* + c_2c_2^*) + 4p \sqrt{p^2-1} \Re(e^{i
\theta} c_1c_2^*) \| \\
&=& \inf \left\| (2p^2 -1) (aa^* + bb^*)
+2p \sqrt{p^2-1} \cos \theta (aa^* - bb^*) 
\right.
\\
&& \qquad\left.
   + 4p \sqrt{p^2-1} \sin \theta \Im(ab^*) \right\|\\
&=& \inf_{p \geq 1, \theta \in \IR} \left\|
(2p^2-1 + 2p \sqrt{p^2-1} \cos \theta) aa^*
\right. \\
&& \qquad + (2p^2-1 - 2p \sqrt{p^2-1} \cos \theta) bb^*
\\
&& \qquad\left.
   + 4p \sqrt{p^2-1} \sin \theta \Im(ab^*) \right\|
\end{eqnarray*}
The claimed formula
follows by taking $r = 2p^2-1 + 2p \sqrt{p^2-1} \cos
\theta$, $s = 2p^2-1 - 2p \sqrt{p^2-1} \cos \theta$ and
$t = 2p \sqrt{p^2-1} \sin \theta$, noting that $rs - t^2 =1$.
We can recover $p$ and $\cos \theta$ from $r,s $ (with $r >0$, $s >
0$, $rs \geq 1$) using
$r + s = 2 (2p^2-1)$, $r-s = 4p \sqrt{p^2-1} \cos \theta$. 
From the sign of $t = \pm \sqrt{rs-1}$ we get $\sin \theta$ and so
$\theta$ modulo $2\pi$.

\begin{remark}
\rm
In \cite{Mag} it is also shown that, for $T$ as in the example above,
$\|T\|_{cb} = \|T\|$. A more general result can be found in
\cite{TimoneyNE}.
\end{remark}

\vfill

\noindent
2000 Mathematics Subject Classification: 47B47, 46L07, 47L25
\vfill

\noindent
School of Mathematics\\
Trinity College\\
Dublin 2\\
Ireland\\[2mm]
Email: \texttt{richardt@maths.tcd.ie}

\end{document}